\documentclass[12pt]{amsart}

\usepackage{amsmath, amsfonts,amsthm,times,graphics}

\usepackage[active]{srcltx}

 \makeatletter
\renewcommand*\subjclass[2][2010]{%
  \def\@subjclass{#2}%
  \@ifundefined{subjclassname@#1}{%
    \ClassWarning{\@classname}{Unknown edition (#1) of Mathematics
      Subject Classification; using '2010'.}%
  }{%
    \@xp\let\@xp\subjclassname\csname subjclassname@#1\endcsname
  }%
}

 \makeatother
\newtheorem{theorem}{Theorem}[section]

\newtheorem{proposition}[theorem]{Proposition}

\theoremstyle{definition}
\newtheorem{definition}[theorem]{Definition}

\newtheorem{remark}[theorem]{Remark}

 \makeatletter
\renewcommand*\subjclass[2][2010]{%
  \def\@subjclass{#2}%
  \@ifundefined{subjclassname@#1}{%
    \ClassWarning{\@classname}{Unknown edition (#1) of Mathematics
      Subject Classification; using '1991'.}%
  }{%
    \@xp\let\@xp\subjclassname\csname subjclassname@#1\endcsname
  }%
}

 \makeatother

\begin{document}
\title[A note on some sub-Gaussian 
 random variables]{A note on some sub-Gaussian 
 random variables}

\author{Romeo Me\v strovi\' c}
\address{Maritime Faculty Kotor, University of Montenegro, 
85330 Kotor, Montenegro} \email{romeo@ac.me}

 
\begin{abstract}
In \cite{m8} the author of this paper continued the research on 
the  complex-valued discrete random variables  $X_l(m,N)$ ($0\le l\le N-1$, 
$1\le M\le N)$  recently introduced and studied  in \cite{ssa}.
Here we extend our results by considering   $X_l(m,N)$
as sub-Gaussian random variables. Our investigation is motivated by the known 
fact that   
 the so-called  Restricted Isometry Property (RIP) 
introduced in \cite{ct1}  holds 
with high probability for any matrix generated by a sub-Gaussian random
variable. Notice that sensing matrices with the RIP play a crucial 
role in   Theory of compressive sensing. 

Our main results concern the proofs of  the lower and upper bound 
 estimates of the expected values 
of the random variables $|X_l(m,N)|$, $|U_l(m,N)|$ and 
$|V_l(m,N)|$, where $U_l(m,N)$ and $U_l(m,N)$ are the real 
and the imaginary part of  $X_l(m,N)$, respectively. 
These estimates are also given in  terms of related sub-Gaussian norm 
$\Vert \cdot\Vert_{\psi_2}$ considered in \cite{ver2}. Moreover,
we prove a refinement of the  mentioned  upper bound 
estimates for the real and the imaginary part of  $X_l(m,N)$.  
   \end{abstract}  
  \maketitle

\section{Introduction and preliminary results}

The recent paper \cite{ssa} by LJ. Stankovi\'c, S. Stankovi\'c
and M. Amin provides a statistical analysis for efficient 
detection of signal components when missing data samples are present
(cf. \cite{sso}, \cite[Section 2]{sid}, \cite{sdsv} and \cite{sdv}).
This analysis is closely related to compressive sensing type problems.
 For more 
information on the development of  compressive sensing (also known as {\it compressed sensing}, 
{\it compressive sampling}, or {\it sparse recovery}), see \cite{do}, \cite{fr}, 
\cite[Chapter 10]{s1} and \cite{sdt}. 
For an excellent survey on this topic with applications
 and related references see  \cite{sr} (also see \cite{op}).
Notice that in the statistical methodology presented in 
\cite{ssa}  a class of complex-valued discrete random variables
(denoted in \cite{m8} as $X_l(m,N)$ with 
 $0\le l\le N-1$ and $1\le M\le N$),  plays  a crucial role.

Following \cite{m8}, the random variable $X_l(m,N)$ can be defined as
follows.

  \begin{definition}(\cite[Definition 1.2]{m8})
Let $N$, $l$ and $m$  be arbitrary nonnegative integers 
such that $0\le l\le N-1$ and $1\le m\le N$.
Let $\Phi(l,N)$ 
be a multiset defined as 
   $$
\Phi(l,N)=\{e^{-j2nl\pi/N}:\, n=1,2,\ldots,N\}.\leqno(1)
   $$
Define the discrete complex-valued random variable $X_l(m,N)=X_l(m)$ as
  \begin{eqnarray*}
 && \mathrm{Prob}\left(X_l(m,N)
 = \sum_{i=1}^me^{-j2n_il\pi/N}\right)\\
(2) &= &\frac{1}{{N\choose m}}\cdot \big|\{\{t_1,t_2,\ldots,t_m\}
\subset\{1,2,\ldots,N\}: 
\sum_{i=1}^me^{-j2t_il\pi/N}=\sum_{i=1}^me^{-j2n_il\pi/N} \big|
\quad\qquad\\
& =&:\frac{q(n_1,n_2,\ldots, n_m)}{{N\choose m}},
  \end{eqnarray*}
  where $\{n_1,n_2,\ldots, n_m\}$ is an arbitrary fixed
 subset of $\{1,2,\ldots,N\}$ such that $1\le n_1<n_2<\cdots <n_m\le N$;
moreover, $q(n_1,n_2,\ldots, n_m)$ is the cardinality of a collection
 of all subsets
$\{t_1,t_2,\ldots,t_m\}$ of the set $\{1,2,\ldots,N\}$ such that 
$\sum_{i=1}^me^{-j2t_il\pi/N}=\sum_{i=1}^me^{-j2n_il\pi/N}$.   
  \end{definition}

Let us recall that by (2) is well defined the random variable
$X_l(m,N)$ taking into account the general additive property of 
probabiblity function $\mathrm{Prob}(\cdot)$ and the fact that there are 
${N\choose m}$ index sets $T\subset \{1,2,\ldots,N \}$ with 
$m$ elements. 

As noticed in \cite[Definition 1.2']{m8}, the random variable $X_l(m,N)$
can be formally  expressed as a sum  
  $$
X_l(m,N)=\sum_{n\in S}e^{-j2nl\pi/N},\leqno(3)
  $$
where the summation ranges over 
any subset  $S$ of size $m$ (the so-called $m$-element subset) 
without replacement from the set $\{1,2,\ldots, N\}$. Notice that the 
number of these subsets $S$ of 
$\{1,2,\ldots, N\}$ is ${N\choose m}$,  and the probability of each  
value of $X_l(m,N)$ is assumed to be equal $1/{N\choose m}$.
   
As usually, throughout our considerations we use the 
term ``multiset'' (often written as ``set'') to mean ``a totality having 
possible multiplicities''; so that two (multi)sets will be counted as equal if 
and only if they have the same elements with identical multiplicities.

Here as always in the sequel, we will denote by
$\Bbb E[X]$ and ${\rm Var}[X]$ the expected value and the variance
of any complex-valued (or real-valued) random variable $X$.
Moreover, for any random variable $X_l(m,N)$ from Definition 1.1
  we shall write
  $$
X_l(m,N)= U_l(m,N)+jV_l(m,N),
 $$
where $U_l(m,N)$ is the {\it real part} 
 and $V_l(m,N)$ is the {\it imaginary part} of $X_l(m,N)$.
Of course,  $U_l(m,N)$ and $V_l(m,N)$ can be considered as 
the real-valued random variables associated with $X_l(m,N)$.
If $l\ge 1$, then it was proved in \cite{ssa} 
(also see \cite[(18) of Theorem 2.4]{m8}) 
that    
  $$
\Bbb E[X_l(m,N)]=\Bbb E[U_l(m,N)]=\Bbb E[U_l(m,N)]=0,\leqno(4) 
   $$  
  Furthermore, it was proved in \cite{ssa} 
(also see \cite[(19) of Theorem 2.4]{m8}) that 
  $$
{\rm Var}[X_l(m,N)]=\Bbb E[|X_l(m,N)|^2]=\frac{m(N-m)}{N-1},\leqno(5)
   $$
whenever $1\le l\le N-1$ and $1\le m\le N$. Moreover, 
if in addition, we suppose that  $N\not= 2l$, then
\cite[(23) of Corollary 2.6]{m8}
  $$
\Bbb E[(U_l(m,N))^2]=\Bbb E[(V_l(m,N))^2]=
\frac{m(N-m)}{2(N-1)}.\leqno(6)
  $$
It was also proved in  \cite[Theorem 2.8]{m8}
that if $l\not= 0$, then for every positive 
integer  $k$ that   is not divisible by 
$N/\gcd(N,l)$  $($$\gcd(N,l)$ denotes  the greatest 
common divisor of $N$  and $l$), the $k$th moment 
$\mu_k:=\Bbb E [X_l(m,N)]$ of the  random variable $X_l(m,N)$  is equal to zero.
In  general case, $\mu_k=\Bbb [X_l(m,N)]$ is a real 
number  \cite[Proposition 2.10]{m8}.

Notice that (1) for $l=0$ implies that
 $$
\Phi(0,N)=\{\underbrace{1,\ldots,1}_N\}.
$$
Moreover, it is obvious that the multiset $\Phi(l,N)$ given by (1) is in fact 
the  set consisting of $N$ (distinct) elements if and only if  $l$ and $N$ 
are relatively prime positive integers (for related discussion, 
see \cite{m5}).

Recall that by using  an Elementary  Number Theory approach
 to some  compressive  sensing problems, different 
classes of random variables $X_l(m,N)$  are considered and 
compared in \cite{m5}.  
Furthermore, in order to establish a probabilistic approach to Welch bound 
on the coherence of a  matrix over the field $\Bbb C$ (or $\Bbb R$),   
a  generalization of the random variable $X_l(m,N)$ is defined and studied
in  \cite{m1}. For more information on the coherence of a  matrix
and related   Welch bound, see \cite[Chapter 5, Theorem 5.7]{fr} 
(also see  \cite{ss}, \cite{s2} and \cite{we}).

Notice also that a {\it Bernoulli probability model}, similar 
to the distribution $\widetilde{X}_l(m,N)$ defined below, was often used in 
the famous paper \cite{crt} by   Cand\`{e}s, Romberg and Tao.
Accordingly, we believe that for  some further probabilistic 
studies  of sparse signal recovery, 
it can be of interest the complex-valued discrete random variable
$\widetilde{X}_l(m,N)$ defined in \cite{m4}. 
Namely, for nonnegative  integers   $N$, $l$ and $m$  such that $0\le l\le N-1$
and $1\le m\le N$, in \cite{m4} it was studied in some sense
 analogous  random variable 
$\widetilde{X}_l(m,N)$ to the random variable ${X}_l(m,N)$, defined as 
a sum
  $$
\widetilde{X}_l(m,N)=\sum_{n=1}^N\exp\left({-\frac{2jn l\pi}{N}}
\right)B_n,
   $$
where  $B_n$ $(n=1,\ldots,N)$ are independent identically distributed  
{\it Bernoulli random variables} ({\it binomial distributions}) taking only
the values 0 and 1 with probability 0 and $m/N$, respectively,
i.e.,
   $$
B_n=\left\{ 
 \begin{array}{ll}
0 & \mathrm{with\,\,probability \,\,} 1-\frac{m}{N}\\
1 & \mathrm{with\,\,probability \,\,} \frac{m}{N}. 
 \end{array}\right.
  $$
Clearly, the range of the random variable 
$\widetilde{X}_l(m,N)$ consists of all possible $2^N-1$ sums of the elements 
of the (multi)set $\{e^{-j2nl\pi/N}:\, n=1,2,\ldots,N\}$.

If $l\ge 1$, then it is proved in \cite[Proposition 2.1]{m4} that
   $$
\Bbb E[\widetilde{X}_l(m,N)]=
\Bbb E[ \widetilde{U}_l(m,N)]=\Bbb E[\widetilde{V}_l(m,N)]=0.\leqno(7) 
   $$
Furthermore, it is proved in  \cite[Proposition 2.1]{m4} that
    $$
{\rm Var}[\widetilde{X}_l(m,N)]=\frac{m(N-m)}{N}.\leqno(8)
   $$
If in addition we suppose that $N\not= 2l$, then \cite[Proposition 2.1]{m4}
      $$
{\rm Var}[\widetilde{U}_l(m,N)]={\rm Var}[\widetilde{V}_l(m,N)]=
\frac{m(N-m)}{2N}.\leqno(9)
   $$

\begin{remark}
From (4) and (7) it follows that for each 
$l\ge 1$ $X_l(m,N)$ and $\widetilde{X}_l(m,N)$
are zero-mean random variables.  
From the expressions (5) and (8) it follows that 
   $$
\frac{{\rm Var}[X_l(m,N)]}{{\rm Var}[\widetilde{X}_l(m,N)]}
=\frac{N}{N-1},\leqno(10)
   $$
i.e.,
 $$
\frac{\sigma [X_l(m,N)]}{\sigma[\widetilde{X}_l(m,N)]}
=\sqrt{\frac{N}{N-1}}.\leqno(11)
   $$
Furthermore, if $N\not= 2l$, then from (6) and (9) of 
\cite[Theorem 2.4]{m8} we find  that 
the proportions (10)  and (11) are also valid 
after replacing  $X_l(m,N)$ by $U_l(m,N)$ (resp. $V_l(m,N)$)
and $\widetilde{X}_l(m,N)$ by $\widetilde{U}_l(m,N)$
(resp. $\widetilde{V}_l(m,N)$).

Notice that in Statistics  the uncorrected sample 
variance or sometimes the 
variance of the sample (observed values) $\{x_1,x_2,\ldots,x_N\}$ 
with the arithmetic mean value $\bar{x}$, is defined as 
       $$
s_N=\frac{1}{N}\sum_{i=1}^N(x_i-\bar{x})^2.\leqno(12) 
    $$
If the biased  sample variance (the second central moment of the sample,
which is a downward-biased estimate of the population variance)
is used to compute an estimate of the population standard deviation,
the result is equal to $s_N$ given by the above formula.

An {\it unbaised estimator} of the variance is given by applying 
{\it Bessel's correction},
using $N-1$ instead of $N$ to yield the {\it unbiased sample variance}, 
denoted by $\bar{s}_N^2$ and defined as
  $$
\bar{s}_N^2=\frac{1}{N-1}\sum_{i=1}^N(x_i-\bar{x})^2.\leqno(13) 
   $$
From (12) and (13) we see that the proportion (10) can be extended
as
   $$
\frac{{\rm Var}[X_l(m,N)]}{{\rm Var}[\widetilde{X}_l(m,N)]}
=\frac{\bar{s}_N^2}{s_N^2}=\frac{N}{N-1}.
   $$
The above proportion suggests the fact that probably in some statistical
sense between the random variables  $X_l(m,N)$ and $\widetilde{X}_l(m,N)$
  there exists a ``connection of type unbiased  sample variance - biased  
sample variance''.  
Moreover,  the values $N/(N-1)$ should be influenced  by the fact that 
  $\widetilde{X}_l(m,N)$ is a sum of $N$ 
independent random variables, while the random variable  $X_l(m,N)$   
is defined on the set 
$\Phi(l,N)$ consisting of $N$ (not necessarily distinct) elements 
that  are ``not independent'' in the sense that their sum 
is equal to zero. 
   \end{remark}

Notice that the random variables  $X_l(m,N)$ and $\widetilde{X}_l(m,N)$
and their real and imaginary parts are bounded random variables.
Therefore, all these  random variables are sub-Gaussian (see Section 2).
In Section 2, we give the assertions  
concerning the  lower and upper bound 
 estimates of the expected values 
of the random variables $|X_l(m,N)|$, $|U_l(m,N)|$ and 
$|V_l(m,N)|$.
These estimates are also given in  terms of related sub-Gaussian norm 
$\Vert \cdot\Vert_{\psi_2}$ considered in \cite{ver2}. Moreover,
we formulate  a refinement of the all mentioned
upper bound  estimates concerning the random variables $|U_l(m,N)|$ and 
$|V_l(m,N)|$. Proofs of all these estimates are given in Section 3.

\section{The main results}

 \begin{theorem} 
Let $N\ge 2$, $l$ and $m$ be nonnegative integers such that 
$0\le l\le N-1$ and  $1\le m\le N$. Then the following 
probability estimates are satisfied:
     \begin{itemize}
\item[(i)] $e^{\frac{N-m}{m(N-1)}}\le 
\Bbb E\left[\exp\left(\frac{\left|X_l(m,N)\right|^2}{m^2}\right)\right]\le e;$
\item[(ii)] $e^{\frac{N-m}{2m(N-1)}}\le 
\Bbb E\left[\exp\left(\frac{(U_l(m,N))^2}{m^2}\right)\right]\le e;$
\item[(iii)] $e^{\frac{N-m}{2m(N-1)}}\le 
\Bbb E\left[\exp\left(\frac{(V_l(m,N))^2}{m^2}\right)\right]\le e$ 
\quad{\rm if}\quad $l\ge 1$.
     \end{itemize}
      \end{theorem}

Notice that the estimates    on the right 
hand side of  (i), (ii) and (iii) of Theorem 2.1  are rough because 
of the fact they are directly obtained by using only the fact that the random variables
$\left|X_l(m,N)\right|$, $\left| U_l(m,N)\right|$ and $\left| V_l(m,N)\right|$
are upper bounded by the constant $m$. Accordingly, 
if $l\ge 1$, then the  equality in each of these inequalities holds if and only 
if $N=1$, i.e., when $X_l(m,N)$, $U_l(m,N)$ nad $V_l(m,N)$ are constant 
random variables  identically equal to one. We believe that for 
non-constant cases, these inequalities  should be 
significantly improved. 

   Theorem 2.1 can be reformulated as follows.
 \begin{theorem} 
Let $N\ge 2$, $l$ and $m$ be nonnegative integers such that 
$0\le l\le N-1$ and  $1\le m\le N$. Then the following 
probability estimates are satisfied:
          \begin{itemize}
\item[(i)] $e^{\frac{m(N-m)}{(N-1)}}\le 
\Bbb E\left[\exp\left(\left|X_l(m,N)\right|^2\right)\right]\le e^{m^2};$
\item[(ii)] $e^{\frac{m(N-m)}{2(N-1)}}\le 
\Bbb E\left[\exp\left((U_l(m,N))^2\right) \right]\le e^{m^2};$
\item[(iii)] $e^{\frac{m(N-m)}{2(N-1)}}\le 
\Bbb E\left[\exp\left((V_l(m,N))^2\right)\right]\le e^{m^2}$ 
\quad{\rm if}\quad $l\ge 1$.
     \end{itemize}
      \end{theorem}

Let us recall that a real-valued random variable $X$ is 
{\it sub-Gaussian} if its distribution is dominated by a normal distribution. 
More precisely, a real-valued random variable $X$ is  sub-Gaussian if 
there holds 
     $$
\mathrm{Prob}(|X|>t)\le \exp\left(1-\frac{t^2}{C^2}\right)
\quad {\rm for \,\, all}\quad t\ge 0,
     $$
where $C>0$ is a real constant that does not depends  on $t$.

 A systematiac introduction into 
sub-Gaussan random variables can be found in \cite[Lemma 5.5 in Subsection 5.2.3 
and Subsection 5.2.5]{ver}; here we briefly mention the basic definitions.
Notice that the {\it  Restricted Isometry Property} (RIP) holds 
with high probability for any matrix generated by a sub-Gaussian random 
variable (see \cite{ct2}, \cite{rv}).
Moreover, a relationship between  the concepts of  coherence and RIP of a 
matrix was established  in \cite{bdf} and \cite{can}. Namely, in these papers 
it is proved that a matrix $A$ with the coherence $\mu(A)$ 
 satisfies the RIP with 
the sparsity order $k\le \frac{1}{\mu(A)}+1$. Therefore, 
it is desirable to give explicit construction 
of matrices with small coherence in compressed sensing.

One of several equivalent ways to define this rigorously is to require 
the {\it Orlicz norm} 
$\Vert X\Vert_{\psi_2}$ defined as 
  $$
\Vert X\Vert_{\psi_2}:=\inf\{ K>0:\, \Bbb 
E\left[\psi_2\left(\frac{|X|}{K}\right)\right]\le 1\}
    $$
to be finite, for the {\it Orlicz function}  $\psi_2(x)=\exp(x^2)-1$. 
The class of sub-Gaussian random variables on a given probability space 
is thus a {\it normed space} endowed with Orlicz norm 
$\Vert \cdot\Vert_{\psi_2}$.  
This definition is in spirit topological in view of the fact that 
the  classical Orlicz norm is used for the definition of  many  
topological vector spaces. For more details on  the Orlicz function and 
related  topological vector spaces,  see 
\cite{mu}. Recall that in Real and Complex Analysis many function spaces are 
endowed with the topology induced by an Orlicz norm (see  \cite[Chapter 7]{mp}
 and \cite{mpl}). 

Obviously, (cf. \cite[Definitions 2.5.6 and Example 2.7.13]{ver2}) 
the  above Orlicz norm $\Vert\cdot \Vert_{\psi_2}$ is exactly the 
{\it sub-Gaussian norm} $\Vert\cdot \Vert_{\psi_2}'$
which is for the sub-Gaussian  real-valued random variable $X$ defined as
   $$
\Vert X \Vert_{\psi_2}'=
\inf\{ K>0:\, \Bbb 
E\left[\exp\left(\frac{X^2}{K^2}\right)\right]\le 2\}.
    $$ 
Accordingly, in the sequel we shall write $\Vert\cdot \Vert_{\psi_2}$
instead of $\Vert\cdot \Vert_{\psi_2}^{'}$.
    
In view of the mentioned facts, a random variable 
$X$ is  sub-Gaussian if and only if
    $$
 \Bbb E\left[\exp\left(\frac{X^2}{\psi}\right)\right]\le 2
   $$
for some real constant $\psi>0$. Hence, any bounded real-valued random 
variable $X$ is sub-Gaussian, and clearly, there holds
   $$
\Vert X \Vert_{\psi_2}\le \frac{1}{\sqrt{\ln 2}}\Vert X \Vert_{\infty}\approx
 1.20112 \Vert X \Vert_{\infty},
  $$
where $\Vert \cdot \Vert_{\infty}$ is the usual  {\it supremum norm}.
Moreover, if $X$ is a centered normal random variable with variance 
$\sigma^2$, then $X$ is sub-Gaussian with 
$\Vert X \Vert_{\psi_2}\le C\sigma$, where $C$ is an absolute constant
\cite[Subsection 5.2.4]{ver}. 
  
Another definition of the sub-Gaussian norm
$\Vert X \Vert_{\psi_2}^{''}$ for the sub-Gaussian  
random variable $X$ was given in \cite[Definition 5.7]{ver} as
   $$
\Vert X \Vert_{\psi_2}^{''}=
\sup_{p\ge 1}\left( p^{-1/2} \, \Bbb (E[|X|^p])^{1/p}\right).
    $$ 
Obviously, there holds
   $$
\Vert X \Vert_{\psi_2}^{''}\le \Vert X \Vert_{\infty}.
  $$

In particular, $X_l(m,N)$, $U_l(m,N)$ and $V_l(m,N)$ are sub-Gaussian
random variables. Clearly, in terms of the sub-Gaussian norm
$\Vert \cdot\Vert_{\psi_2}$ Theorem 2.2 can be reformulated as follows.

\begin{proposition}
Let $N\ge 1$, $l$ and $m$ be nonnegative integers such that 
$0\le l\le N-1$ and  $1\le m\le N$. 
Then $\left|X_l(m,N)\right|$, $U_l(m,N)$ and
 $V_l(m,N)$ are  sub-Gaussian random variables. Moreover, there holds
    \begin{itemize}
\item[(i)]
$\sqrt{\frac{m(N-m)}{(N-1)\ln 2}}\le \Vert |X_l(m,N)|\Vert_{\psi_2}\le 
\frac{m}{\sqrt{\ln 2}};$

\item[(ii)] $\sqrt{\frac{m(N-m)}{2(N-1)\ln 2}}\le \Vert U_l(m,N)\Vert_{\psi_2}\le 
\frac{m}{\sqrt{\ln 2}};$
\item[(iii)]
$\sqrt{\frac{m(N-m)}{2(N-1)\ln 2}}\le \Vert V_l(m,N)\Vert_{\psi_2}\le 
\frac{m}{\sqrt{\ln 2}}\quad{\rm if}\quad l\ge 1$.
    \end{itemize}
      \end{proposition}

The upper bound $m/\sqrt{\ln 2}$ on the right hand side of the estimates 
(ii) and (iii)
of Proposition 2.3 can be improved for a large class of random variables 
$U_l(m,N)$ and $V_l(m,N)$. This is given by the following result.

  \begin{proposition}
Let $N\ge 2$, $l$ and $m$ be positive  integers such that 
$1\le l\le N-1$ and  $1\le m\le N$. If $N$ and $l$ are relatively prime,
then 
  $$
 \Vert U_l(m,N)\Vert_{\psi_2}\le \frac{\sin\frac{m\pi }{N}}
{\sqrt{\ln 2} \sin\frac{\pi }{N}}\leqno(14)
    $$
and
     $$
\Vert V_l(m,N)\Vert_{\psi_2}\le \left\{ 
 \begin{array}{ll}
  \frac{\sin\frac{m\pi}{N}\sin\frac{(2\lfloor N/4\rfloor+1)\pi}{N}}
{\sqrt{\ln 2}\sin\frac{\pi }{N}}&  if\,\, m\,\,  is\,\, even\\
\frac{\sin\frac{m\pi}{N}\sin\frac{2\lfloor (N+1)/4\rfloor\pi}{N} }
{\sqrt{\ln 2}\sin\frac{\pi }{N}}&  if\,\, m\,\,  is\,\, odd.
\end{array}\right. \leqno(15)
    $$
    \end{proposition}
 \begin{remark}
Notice that  if $m\sim cN$ for some constant $c$ with $0<c\le 1/2$ and all 
sufficiently large values of $N$, then $\sin (\pi/N)\approx \pi/N$ and 
thus,  the upper bound
on the right hand side of estimates (14) and (15) is 
$$
\sim N\sin (c\pi)/(\pi\sqrt{\ln 2})=0.382329 N\sin (c\pi).
$$
On the other hand, from (ii) and (iii) of Proposition 2.3  we see that  for 
such a value $m$,  the lower bound on the left hand side of 
the estimates (ii) and (iii)
is 
$$
\sim \sqrt{c(1-c)N/(2\ln 2)}.
 $$
For example, if $m\sim N/2$ (i.e., $c=1/2$), then these upper and lower bounds 
are approximately equal to  $0.382329 N$ and $0.424661\sqrt{N}$, 
respectively.
\end{remark}

From the estimates (14), (15) and proof of Proposition 2.4,   taking into
account that $|X_l(m,N)|=\sqrt{(U_l(m,N))^2+ (V_l(m,N))^2}$,
it follows immediately the following result.

  \begin{proposition}
Let $N\ge 2$, $l$ and $m$ be positive  integers such that 
$1\le l\le N-1$ and  $1\le m\le N$. If $N$ and $l$ are relatively prime,
then 
  $$
 \Vert X_l(m,N)\Vert_{\psi_2}\le \frac{\sqrt{2}\sin\frac{m\pi }{N}}
{\sqrt{\ln 2} \sin\frac{\pi }{N}}.
    $$
  \end{proposition}

\section{Proofs of the results}

\begin{proof}[Proof of Theorem $2.1$]
First notice that for $l=0$ and all $m$ with $1\le m\le N$, 
$\left|X_0(m,N)\right|$, $U_0(m,N)$ and $V_0(m,N)$ are
 constant random variables  with 
 $$
{\rm Prob}\left(\left|X_0(m1,N)\right|=m\right)=
{\rm Prob}\left(U_0(m,N)=m\right)=
{\rm Prob}\left(V_0(m,N)=0\right)=1.
 $$
Therefore, both double inequalities (i) and  (ii)  are satisfied.   

Now suppose that $1\le l\le N-1$.
Since the random variables  $|X_l(m,N)|^2$, $(U_l(m,N))^2$
and $(V_l(m,N))^2$ are obviously bounded below by the constant $m^2$,
the inequalities on the right hand side of (i), (ii) and (iii)
are trivially satisfied.

  Notice that  
  \begin{equation*}\begin{split}
(16)&\Bbb E\left[\exp\left(\frac{|X_l(m,N)|^2}{m^2}\right)\right]\\
&=\frac{1}{{N\choose m}} \Big( \sum_{\{i_1,i_2,\ldots,i_m\} 
\subset \{1,2,\ldots,N\}}\exp\frac{(w_{i_1}+w_{i_2}+\cdots+
w_{i_m})(\overline{w_{i_1}+w_{i_2}+\cdots+
w_{i_m}})}{m^2}\Big),
   \end{split}\end{equation*}
where the summation ranges over all 
${N\choose m}$ subsets $\{i_1,i_2,\ldots,i_m\}$
of $\{1,2,\ldots,N\}$ with $1\le i_1<i_2<\cdots <i_m\le N$.
Notice that 
   $$
A_{\{i_1,i_2,\ldots,i_m\}}:=\exp\big( (w_{i_1}+w_{i_2}+\cdots+
w_{i_m})(\overline{w_{i_1}+w_{i_2}+\cdots+
w_{i_m}})\big)
  $$ 
are positive real numbers for  each subset $\{i_1,i_2,\ldots,i_m\}$
of $\{1,2,\ldots,N\}$ with $1\le i_1<i_2<\cdots <i_m\le N$.
Then applying to these numbers the classical {\it arithmetic-geometric mean 
inequality} $(\sum_{k=1}^na_k)/n\ge \sqrt[n]{\prod_{k=1}^na_k}$
($n\in \Bbb N$, $a_1,\ldots, a_n\in \Bbb R^{+}$), and using the expression 
(16), we find that the right hand side
of this expression  is
   \begin{equation*}\begin{split}
\ge & \sqrt[{N\choose m}]{\exp\left(\frac{1}{m^2}\sum_{\{i_1,i_2,\ldots,i_m\} 
\subset \{1,2,\ldots,N\}} (w_{i_1}+w_{i_2}+\cdots+
w_{i_m})(\overline{w_{i_1}+w_{i_2}+\cdots+
w_{i_m}})\right)}\\
&= \sqrt[{N\choose m}]{\exp \Big( \frac{1}{m^2}{N\choose m}\Bbb 
E\left[|X_l(m,N)|^2\right]\Big)}=\exp\left(\frac{N-m}{m(N-1)}\right). 
   \end{split}\end{equation*}
This proves the left hand side of the inequality (i) 
of Theorem 2.1.

Finally, notice that the left hand sides of  inequalities (ii) 
and (iii) of Theorem 2.1 can be proved in the same manner as that of (i),
 using in the final step the first and the second equality
of the  expression  (6), respectively.
Hence, these proofs can be omitted, and proof of Theorem 2.1 
is completed.     
   \end{proof}
\begin{proof}[Proof of Theorem $2.2$]
Proof of Theorem 2.2 is completely similar to those of Theorem 2.1
and hence, may be omitted.
  \end{proof}

\begin{proof}[Proof of Proposition $2.3$]
The first assertion is an immediate consequence of inequalities on 
the right hand sides of (i), (ii) and (iii) of Theorem 2.1.
 The  inequalities on the right hand side of (i), (ii) and (iii) 
are also  immediate consequences of  the inequalities on 
the right hand sides of (i), (ii) and (iii) of Theorem 2.1, respectively.
 Finally, the  inequalities on the left hand side of (i), (ii) and (iii) 
 can be proved in the same manner as those of (i)
of Theorem 2.1.
\end{proof}

   \begin{proof}[Proof of Proposition $2.4$]
Since by the assumption, $N$ and $l$ are relatively prime positive integers,
then the multiset $\Phi(l,N)$ defined by (1) consists  of
$N$ distinct elements, and it can be written as    
  $$
\Phi(l,N)=\{1,w,w^2,\ldots,w^{N-1}\},\leqno(17)
  $$
where $w=\exp\left(2\pi j/N\right)$ is the primitive $N$th root of 
unity. Then the ranges (the sets of all values) of the random variables 
$U_l(m,N)$ and $V_l(m,N)$ are respectively given by
    \begin{equation*}\begin{split}
(18)  &{\mathcal R}(U_l(m,N))\\
 &=\left\{\cos\frac{2k_1\pi}{N}+
\cos\frac{2k_2\pi}{N}+\cdots +\cos\frac{2k_m\pi}{N}:
0\le k_1<k_2<\cdots<k_m\le N-1\right\}
     \end{split}\end{equation*}
and 
   \begin{equation*}\begin{split}
(19) & {\mathcal R}(V_l(m,N))\\
\qquad &=\left\{\sin\frac{2k_1\pi}{N}+
\sin\frac{2k_2\pi}{N}+\cdots +\sin\frac{2k_m\pi}{N}:
0\le k_1<k_2<\cdots<k_m\le N-1\right\}.
     \end{split}\end{equation*} 
In the whole  proof $M_1$ and $M_2$ will always denote 
the maximal value and the minimal value of considered 
 random variable $U_l(m,N)$ or  $V_l(m,N)$, respectively. 
In order to obtain the upper bounds for  $\Vert U_l(m,N)\Vert_{\infty}$
and $\Vert V_l(m,N)\Vert_{\infty}$,
in view of the antisymmetric property of random variables $U_l(m,N)$
and $V_l(m,N)$  given in  \cite[Proposition 2.1]{m8}, without loss of generality, 
in the whole proof we can suppose that 
$m\le \lfloor N/2\rfloor$ ($\lfloor x\rfloor$ denotes the greatest integer 
not exceeding a real number $x$).

{\it Proof of the inequality $(14)$}.
As noticed in Section 2,  every bounded 
random variable $X$ is sub-Gaussian, and  there holds 
     $$
\Vert X\Vert_{\psi_2}\le \frac{1}{\sqrt{\ln 2}}\Vert X\Vert_{\infty},\leqno(20)
    $$
where $\Vert \cdot \Vert_{\infty}$ is the usual  supremum norm.

We will consider the cases when a positive integer  $m$ is odd and 
when $m$ is even.

{\it The first case}:  $m$ is an odd positive integer. Put 
$m=2s+1$ with integer $s\ge 0$. If $s=0$ then $m=1$, and hence, 
    $$
{\mathcal R}(U_l(1,N))=\left\{1,\cos\frac{2\pi}{N},\ldots,
\cos\frac{2(N-1)\pi}{N}\right\}.
   $$
Therefore, $\Vert U_l(m,N)\Vert_{\infty}\le 1$, which together with 
(20) yields
   $$
\Vert X\Vert_{\psi_2}\le \frac{1}{\sqrt{\ln 2}}.
    $$
Notice that the above inequality coincides with (14) for $m=1$.
    
Now suppose that $s\ge 1$, i.e., $m\ge 3$.
Since by the above assumption, $m\le \lfloor N/2\rfloor$,
it follows that  $s\le \lfloor N/2\rfloor/2-1\le N/4-1$, 
and hence, we have  
      $$
 \cos\frac{2k\pi}{N}>0 \quad {\rm for\,\,all}\quad k=1,2,\ldots,s.
\leqno(21)   
      $$
Since the function 
$f(x):=\cos x$ is decreasing on the segment $[0,\pi]$ and since 
$\cos x=\cos (2\pi-x)$, in view of (18) and (21), we conclude that
the random variable $U_l(m.N)$ attains its maximal value equals to 
  $$
M_1=1+\sum_{k=1}^s \cos\frac{2k\pi}{N}+
\sum_{k=1}^s \cos\frac{2(N-k)\pi}{N}=
1+2\sum_{k=1}^s \cos\frac{2k\pi}{N}.\leqno(22)
  $$ 
Since $\cos\frac{2k\pi}{N}=\Re \left(\exp\left(2k\pi j/N\right)\right)=
\Re (w^k)$, from (22) we obtain
  \begin{equation*}\begin{split}
M_1&=1+ 2\sum_{k=1}^s \Re (w^k)=1+ 2\Re\left(\sum_{k=1}^s w^k\right)\\
&=1+2\Re\left(\frac{w-w^{s+1}}{1-w}\right)=
1+2\Re\left(\frac{w-w^{s+1}}{1-w}\cdot\frac{1-\bar{w}}{1-\bar{w}}\right)\\
 &=1+2\cdot\frac{\Re\left(w-1-w^{s+1}+w^s\right)}{1-2\Re(w)+|w|^2}=
1+2\cdot\frac{\cos\frac{2\pi}{N}-1- \cos\frac{2(s+1)\pi}{N}
+\cos\frac{2s\pi}{N}}{2-2\cos\frac{2\pi}{N}}\\
(23)\quad&=\frac{\cos\frac{2s\pi}{N}-
\cos\frac{2(s+1)\pi}{N}}{1-\cos\frac{2\pi}{N}}\\
&({\rm by\,\, using\,\, the\,\, identities}\,\, 
\cos\alpha-\cos\beta=2\sin\frac{\alpha+\beta}{2}\sin
\frac{\beta-\alpha}{2}\quad{\rm and}\\
& 1-\cos 2\alpha=2\sin^2\alpha)\\
&=\frac{2\sin\frac{(2s+1)\pi}{N}\sin\frac{\pi}{N}}
{2\sin^2\frac{\pi}{N}}=\frac{\sin\frac{(2s+1)\pi}{N}}{\sin\frac{\pi}{N}}
=\frac{\sin\frac{m\pi}{N}}{\sin\frac{\pi}{N}}.
     \end{split}\end{equation*}
In order to determine the minimal value $M_2$ of the random variable 
$U_l(m,N)$,  we will consider the following two subcases:

{\it The first subcase}:  $N$ is an even positive integer.
Take $N=2n$ with $n\in\Bbb N$. Then by using the same argument 
applied for determining the above maximal value $M_1$ of $U_l(m,N)$, (22)  and 
(23), we obtain 
   \begin{equation*}\begin{split}
M_2&=\cos\frac{2n\pi}{2n}+\sum_{t=n-s}^{n-1} \cos\frac{2t\pi}{2n}+
\sum_{t=n+1}^{n+s} \cos\frac{2t\pi}{2n}\\
(24)\qquad\qquad&=-1+\sum_{k=1}^{s} \cos\frac{2(n-k)\pi}{2n}+
\sum_{k=1}^{s} \cos\frac{2(n+k)\pi}{2n}\\
&=-1-\sum_{k=1}^{s} \cos\frac{2k\pi}{2n}-
\sum_{k=1}^{s} \cos\frac{2k\pi}{2n}\quad ({\rm the\,\, change}\quad 2n=N)
\qquad\qquad\\
&=-M_1=-\frac{\sin\frac{m\pi}{N}}{\sin\frac{\pi}{N}}.
    \end{split}\end{equation*}

{\it The second subcase}:  $N$ is an odd positive integer.
Take $N=2n+1$ with $n\in\Bbb N$. Then similarly as above, we find that 
   \begin{equation*}\begin{split}
M_2=&\cos\frac{2(n-s)\pi}{2n+1}+\sum_{t=n-s+1}^{n} \cos\frac{2t\pi}{2n+1}
+\sum_{t=n+1}^{n+s} \cos\frac{2t\pi}{2n+1}\\
&=-\cos\left(\pi- \frac{2(n-s)\pi}{2n+1}\right)\\
&  -\sum_{t=n-s+1}^{n} \cos\left(\pi-\frac{2t\pi}{2n+1}\right)-
\sum_{t=n+1}^{n+s} \cos\left(\frac{2t\pi}{2n+1}-\pi\right)\\
(25)\qquad\qquad\qquad=& -\cos \frac{(2s+1)\pi}{2n+1}
  -\sum_{t=n-s+1}^{n} \cos\frac{(2n+1-2t)\pi}{2n+1}\qquad\qquad\qquad\\
&-\sum_{t=n+1}^{n+s} \cos\frac{(2t-2n-1)\pi}{2n+1}\\
=& -\cos \frac{(2s+1))\pi}{2n+1}
  -\sum_{k=1}^{s} \cos\frac{(2k-1)\pi}{2n+1}-
\sum_{k=1}^{s} \cos\frac{(2k-1)\pi}{2n+1}\\
=& -\cos \frac{(2s+1)\pi}{2n+1} -2\sum_{k=1}^{s} \cos\frac{(2k-1)\pi}{2n+1}.
   \end{split}\end{equation*}
If we put $\xi=\exp\left(j\pi/(2n+1)\right)$, then 
   $\cos \frac{t\pi}{2n+1}=\Re (\xi^t)$ for each $t\in \Bbb N$, 
and hence, from (25) we get
     \begin{equation*}\begin{split}
M_2=& -\cos \frac{(2s+1)\pi}{2n+1} -2\sum_{k=1}^{s}\Re(\xi^{2k-1})=
-\cos \frac{(2s+1)\pi}{2n+1} -2\Re\left(\sum_{k=1}^{s}\xi^{2k-1}\right)\\
=&-\cos \frac{(2s+1)\pi}{2n+1}-2\Re
\left(\frac{\xi-\xi^{2s+1}}{1-\xi^2} \right)\\
=&-\cos \frac{(2s+1)\pi}{2n+1}-2\Re
\left(\frac{\xi-\xi^{2s+1}}{1-\xi^2}\cdot\frac{1-\bar{\xi}^2}
{1-\bar{\xi}^2} \right)\\
=&-\cos \frac{(2s+1)\pi}{2n+1}-2\Re
\left(\frac{\xi-\bar{\xi}-\xi^{2s+1}+\xi^{2s-1}}{1-2\Re(\xi^2)+
|\xi|^4}\right)\\
=&-\cos \frac{(2s+1)\pi}{2n+1}-\frac{2\Re(\xi^{2s-1}-\xi^{2s+1})}{2-2
\Re(\xi^2)}\\
(26)\quad=&-\cos \frac{(2s+1)\pi}{2n+1}-\frac{\cos\frac{(2s+1)\pi}{2n+1} -
\cos\frac{(2s-1)\pi}{2n+1}}{1-\cos\frac{2\pi}{2n+1}}\\
    &({\rm by\,\, using\,\, the\,\, identity}\,\, 
     \cos\alpha-\cos\beta=2\sin\frac{\alpha+\beta}{2}\sin
      \frac{\beta-\alpha}{2}\quad{\rm and}\\
& 1-\cos 2\alpha=2\sin^2\alpha)\\
=&-\cos \frac{(2s+1)\pi}{2n+1}-\frac{2\sin\frac{2s\pi}{2n+1}
\sin\frac{\pi}{2n+1}}
{2\sin^2\frac{\pi}{2n+1}}=-\cos \frac{(2s+1)\pi}{2n+1}-
\frac{\sin\frac{2s\pi}{2n+1}}{\sin\frac{\pi}{2n+1}}\\
=& -\frac{\sin\frac{\pi}{2n+1}\cos\frac{(2s+1)\pi}{2n+1}+
\sin\frac{2s\pi}{2n+1}}{\sin\frac{\pi}{2n+1}}\\
&({\rm by\,\, using\,\, the\,\, identity}\,\, 
   \sin(\alpha -\beta)=\sin\alpha\cos\beta-\cos\alpha\sin\beta)\\
=& -\frac{\sin\frac{\pi}{2n+1}\cos\frac{(2s+1)\pi}{2n+1}+
\sin\frac{(2s+1)\pi}{2n+1}\cos\frac{\pi}{2n+1}-  
\cos\frac{(2s+1)\pi}{2n+1}\sin\frac{\pi}{2n+1}}{\sin\frac{\pi}{2n+1}}\\
=& -\frac{\sin\frac{(2s+1)\pi}{2n+1}\cos\frac{\pi}{2n+1}}
{\sin\frac{\pi}{2n+1}}=-\frac{\sin\frac{(2s+1)\pi}{N}\cos\frac{\pi}{N}}
{\sin\frac{\pi}{N}}=-\frac{\sin\frac{m\pi}{N}\cos\frac{\pi}{N}}
{\sin\frac{\pi}{N}}.
  \end{split}\end{equation*}
 From (23), (24) and (26) we see that 
$|M_2|\le M_1$ for every odd integer 
$m\ge 3$, and hence for such  a $m$ we have  
   $$
\Vert U_l(m,N) \Vert_{\infty} =\max\left\{M_1, |M_2| \right\}
=\frac{\sin\frac{m\pi}{N}}{\sqrt{\ln 2}\sin\frac{\pi}{N}}.\leqno(27)
   $$
From (20) and (27) we immediately obtain
  $$
\Vert U_l(m,N)\Vert_{\psi_2}\le 
\frac{\sin\frac{m\pi}{N}}{\sqrt{\ln 2}\sin\frac{\pi}{N}},\leqno(28)
  $$
as asserted.

{\it The second case}:  $m$ is an even positive integer.
Take $m=2s$ with integer  $s\ge 1$. Then by using the same argument 
applied in the first case, similarly as in the first case, 
we find that the random variable $U_l(m,N)$ attains its maximal value 
equals to 
    \begin{equation*}\begin{split}
M_1=& 1+\cos\frac{2s\pi }{N}+\sum_{k=1}^{s-1} \cos\frac{2k\pi}{N}
+\sum_{k=1}^{s-1} \cos\frac{2(N-k)\pi}{N}\\
=&1+\cos\frac{2s\pi }{N}+2\sum_{k=1}^{s-1} \cos\frac{2k\pi}{N}\\
(29)\qquad\qquad\qquad\qquad
=&\frac{\sin\frac{\pi}{N}+\sin\frac{\pi}{N}\cos\frac{2s\pi}{N}+
\sin\frac{(2s-1)\pi}{N}-\sin\frac{\pi}{N}}{\sin\frac{\pi}{N}}\qquad\qquad
\qquad\qquad\\
=& \frac{\sin\frac{\pi}{N}\cos\frac{2s\pi}{N}+
 \sin\frac{2s\pi}{N}\cos\frac{\pi}{N}-\sin\frac{\pi}{N}
\cos\frac{2s\pi}{N}}{\sin\frac{\pi}{N}}\\
=& \frac{\sin\frac{2s\pi}{N}\cos\frac{\pi}{N}}{\sin\frac{\pi}{N}}=
\frac{\sin\frac{m\pi}{N}\cos\frac{\pi}{N}}{\sin\frac{\pi}{N}}.
   \end{split}\end{equation*} 
If $N=2n$ $(n\in\Bbb N)$ is an even positive integer, then  proceeding in the 
same manner as in the above first subcase (see (24)), 
we obtain that the minimal
value of the random variable $U_l(m,N)$ is equal to
     \begin{equation*}\begin{split}
M_2=& \cos\frac{2n \pi }{2n}+2\sum_{k=n-s+1}^{n-1}\cos\frac{2k\pi}{2n}+
\cos\frac{2(n-s)\pi}{2n}\\
=&-1+\frac{2\sin\frac{(s-1)\pi}{2n}\cos\frac{(2n-s)\pi}{2n}}
{\sin\frac{\pi}{2n}}+\cos\left(\pi-\frac{s\pi}{n}\right)\\
=&-1-\frac{2\sin\frac{(s-1)\pi}{2n}\cos\frac{s\pi}{2n}}
{\sin\frac{\pi}{2n}}-\cos\frac{s\pi}{n}\\
(30)\qquad\qquad\qquad\qquad=&\frac{-\sin\frac{\pi}{2n}- 
2\sin\frac{(s-1)\pi}{2n}\cos
\frac{s\pi}{2n}- \sin\frac{\pi}{2n}\cos\frac{ s\pi}{n}}
{\sin\frac{\pi}{2n}}\qquad\qquad\qquad\\
=&\frac{-\sin\frac{\pi}{2n}-\sin\frac{(2s-1)\pi}{2n}
+\sin\frac{\pi}{2n}-\sin\frac{\pi}{2n}\cos\frac{ s\pi}{n}}
{\sin\frac{\pi}{2n}}\\
=&\frac{-\sin\frac{ 2s\pi}{2n}\cos\frac{\pi}{2n}
+\sin\frac{\pi}{2n}\cos\frac{ s\pi}{n}- 
\sin\frac{\pi}{2n}\cos\frac{ s\pi}{n}}{\sin\frac{\pi}{2n}}\\
=& -\frac{\sin\frac{ 2s\pi}{2n}\cos\frac{\pi}{2n}}{\sin\frac{\pi}{2n}}=
-\frac{\sin\frac{ m\pi}{N} \cos\frac{\pi}{N}}{\sin\frac{\pi}{N}}.
    \end{split}\end{equation*}

 If $N=2n+1$ ($n\in\Bbb N$) is an odd positive integer,  
then similarly as in the previous cases, we obtain that the minimal
value of the random variable $U_l(m,N)$ is equal to
  \begin{equation*}\begin{split}
M_2=& \sum_{k=n-s+1}^{n+s}\cos\frac{2k\pi}{N}=
\frac{\sin\frac{2s\pi}{2n+1}\cos\frac{(2n+1)\pi}{2n+1}}
{\sin\frac{\pi}{2n+1}}\\
(31)\qquad\qquad\qquad\qquad =&-\frac{\sin\frac{m\pi}{N}}{\sin\frac{\pi}{N}}.
\qquad\qquad\qquad\qquad\qquad\qquad\qquad\qquad\qquad
\qquad\qquad\qquad\qquad\qquad\qquad\qquad
  \end{split}\end{equation*}
From (29), (30) and (31) we see that for each 
even integer $m\ge 2$,
  $$
\Vert X_l(m,N)\Vert_{\infty}=\max\{ M_1,|M_2|\} \le  
\frac{\sin\frac{m\pi}{N}}{\sin\frac{\pi}{N}},
  $$ 
which in view of the inequality (20) yields
$$
\Vert U_l(m,N)\Vert_{\psi_2}\le 
\frac{\sin\frac{m\pi}{N}}{\sqrt{\ln 2}\sin\frac{\pi}{N}}.
  $$
Therefore, proof of the inequality (14) is completed.\\

{\it Proof of the inequality $(15)$}.
In order to prove the inequality (15), we proceed similarly as in the case 
of $U_l(m,N)$. Since 
$\sin\frac{2k\pi}{N}=\Im \left(\exp\left(2k\pi j/N\right)\right)=
\Im (w^k)$,  proceeding by the analogus way as in (23) (replacing 
$\Re(\cdot)$ by $\Im(\cdot)$), 
we obtain the following known identity:
   $$
\sum_{k=t}^{t+q}\sin\frac{2k\pi}{N}=
\frac{\sin\frac{(q+1)\pi}{N}\sin\frac{(2t+q)\pi}{N} }{\sin\frac{\pi}{N}},
\leqno(32)
  $$
where $t\ge 1$ and $q\ge 0$ are  nonnegative integers.
Using the identity (32) and considering the cases when $m$ is odd 
and $m$ is even both divided into the following fourt subcases:
$N\equiv 0(\bmod (\,4)$, $N\equiv 1(\bmod (\,4)$, $N\equiv 2(\bmod (\,4)$  
and $N\equiv 3(\bmod (\,4)$, we can arrive at the estimate given by (15)
by considering the  following four cases.

{\it The first case}: $m$ is an even positive integer and $N\equiv 1\pmod{4}$.
Put $m=2s$ and  $N=4n+1$ for some integers $s\ge 1$  and $n\ge 1$.
Then it is  easy to see that 
  $$
M_1=\sum_{k=n-s+1}^{n+s}\sin\frac{ 2k\pi}{4n+1},
  $$  
which by using the identity (32) 
immediately yields
  $$
M_1=\frac{\sin\frac{2s\pi}{4n+1}\sin\frac{(2n+1)\pi}{4n+1}}
{\sin\frac{\pi}{4n+1}}=\frac{\sin\frac{m\pi}{N}
\sin\frac{(2\lfloor N/4\rfloor +  1)\pi}{N}}
{\sin\frac{\pi}{N}}.\leqno(33)
  $$ 
Similarly, we have
 $$
M_2=\sum_{k=3n-s+1}^{3n+s}\sin\frac{ 2k\pi}{4n+1},
  $$  
whence by using the identity (32) 
 it follows that
  $$
M_2=\frac{\sin\frac{2s\pi}{4n+1}\sin\frac{(6n+1)\pi}{4n+1}}
{\sin\frac{\pi}{4n+1}}=-\frac{\sin\frac{2s\pi}{4n+1}\sin\frac{(2n+1)\pi}{4n+1}}
{\sin\frac{\pi}{4n+1}}   =-\frac{\sin\frac{m\pi}{N}
\sin\frac{(2\lfloor N/4\rfloor+1)\pi}{N}}
{\sin\frac{\pi}{N}}.\leqno(34)
  $$ 
From (33) and (34) we immediately obtain
     $$
\Vert V_l(m,N) \Vert_{\infty} =\max\left\{M_1, |M_2| \right\}=
\frac{\sin\frac{m\pi}{N}
\sin\frac{(2\lfloor N/4\rfloor+1)\pi}{N}}
{\sin\frac{\pi}{N}}.\leqno(35)
  $$

{\it The second case}: $m$ is an even positive integer and $N\equiv 3\pmod{4}$.
Put $m=2s$ and  $N=4n+3$ for some integers $s\ge 1$  and $n\ge 1$.
Then as in the first case, it is  easy to see that 
  $$
M_1=\sum_{k=n-s+1}^{n+s}\sin\frac{ 2k\pi}{4n+1},
  $$  
which by using the identity (32) 
immediately yields
  $$
M_1=\frac{\sin\frac{2s\pi}{4n+3}\sin\frac{(2n+1)\pi}{4n+3}}
{\sin\frac{\pi}{4n+3}}=\frac{\sin\frac{m\pi}{N}
\sin\frac{(2\lfloor N/4\rfloor +  1)\pi}{N}}
{\sin\frac{\pi}{N}}.\leqno(36)
  $$ 
Similarly, we have
 $$
M_2=\sum_{k=3n-s+3}^{3n+s+2}\sin\frac{ 2k\pi}{4n+3},
  $$  
whence by using the identity (32), 
it follows that
  $$
M_2=\frac{\sin\frac{2s\pi}{4n+3}\sin\frac{(6n+5)\pi}{4n+3}}
{\sin\frac{\pi}{4n+3}}=-\frac{\sin\frac{2s\pi}{4n+3}\sin\frac{(2n+2)\pi}{4n+3}}
{\sin\frac{\pi}{4n+3}}   =-\frac{\sin\frac{m\pi}{N}
\sin\frac{(2\lfloor N/4\rfloor +1)\pi}{N}}
{\sin\frac{\pi}{N}}.\leqno(37)
  $$ 
The equalities (36) and (37) imply that
     $$
\Vert V_l(m,N) \Vert_{\infty} =\max\left\{M_1, |M_2| \right\}=
\frac{\sin\frac{m\pi}{N}
\sin\frac{(2\lfloor N/4\rfloor+1)\pi}{N}}
{\sin\frac{\pi}{N}}.\leqno(38)
  $$

{\it The third case}: $m$ and $N$ are  even positive integers.
Put $m=2s$ and  $N=2n$ for some integers $s\ge 1$  and $n\ge 1$.
Then it  is easy to check that 
  $$
M_1=\sum_{k=\lfloor n/2\rfloor -s+1}^{\lfloor n/2\rfloor +s}
\sin\frac{ 2k\pi}{2n},
  $$  
which by applying the identity (32) 
and some basic trigonometric identities to  both cases 
$N\equiv 0 \pmod{4}$ and $N\equiv 2 \pmod{4}$,  
we immediately obtain 
  $$
M_1=\frac{\sin\frac{m\pi}{N}
\sin\frac{(2\lfloor N/4\rfloor +  1)\pi}{N}}
{\sin\frac{\pi}{N}}.\leqno(39)
  $$ 
Similarly, we find that
  $$
M_2=\sum_{k=\lfloor 3n/2\rfloor-s+1}^{\lfloor 3n/2\rfloor+s}
\sin\frac{ 2k\pi}{2n},
  $$  
whence by applying the identity (32) 
and some basic trigonometric identities 
we get
  $$
M_2=-\frac{\sin\frac{2s\pi}{2n}\sin\frac{(3n+1)\pi}{2n}}
{\sin\frac{\pi}{2n}}=\frac{\sin\frac{m\pi}{N}
\sin\frac{(2\lfloor N/4\rfloor +  1)\pi}{N}}
{\sin\frac{\pi}{N}}.\leqno(40)
  $$
The equalities (39) and (40) imply that
     $$
\Vert V_l(m,N) \Vert_{\infty} =\max\left\{M_1, |M_2| \right\}=
\frac{\sin\frac{m\pi}{N}
\sin\frac{(2\lfloor N/4\rfloor+1)\pi}{N}}
{\sin\frac{\pi}{N}}.\leqno(41)
  $$

{\it The fourth case}: $m\ge 1$ is an odd  positive integer.
If we take $m=2s+1$ with some integer $s\ge 0$, then 
by considering the all four subcases $N\pmod{4}$, we can easily arrive 
to the equality
    $$
M_1=\sum_{k=\lfloor (N+1)/4\rfloor -s}^{\lfloor (N+1)/4\rfloor +s}
\sin\frac{ 2k\pi}{N},
  $$  
which by applying the identity (32) and  
 some basic trigonometric identities, 
immediately yields
  $$
M_1=\frac{\sin\frac{m\pi}{N}\sin\frac{2\lfloor (N+1)/4\rfloor\pi}{N}}
{\sin\frac{\pi}{N}}.\leqno(42)
  $$
Similarly, we find that
  $$
M_2=\sum_{k=\lfloor (3N+1)/4\rfloor -s}^{\lfloor (3N+1)/4\rfloor +s}
\sin\frac{ 2k\pi}{N},
  $$  
which by applying the identity (32) and 
 some basic trigonometric identities,    
immediately gives
  $$
M_2=\frac{\sin\frac{m\pi}{N}
\sin\frac{2\lfloor (3N+1)/4\rfloor\pi}{N}}
{\sin\frac{\pi}{N}}.\leqno(43)
  $$
If $N$ is even, then $2\lfloor (3N+1)/4\rfloor-2\lfloor (N+1)/4\rfloor=N$,
and thus, from (42) and (43) we have that $M_2=-M_1$.
If $N$ is odd, then $2\lfloor (3N+1)/4\rfloor +2\lfloor (N+1)/4\rfloor=2N$,
and so, from (42) and (43) we also have that $M_2=-M_1$.  
Therefore, for each $N\ge 2$ there holds 
      $$
\Vert V_l(m,N) \Vert_{\infty} =\max\left\{M_1, |M_2| \right\}=M_1=
\frac{\sin\frac{m\pi}{N}
\sin\frac{2\lfloor (N+1)/4\rfloor\pi}{N}}
{\sin\frac{\pi}{N}}.\leqno(44)
  $$
Finally, (20) and the equalities (35), (38), (41) and (44) immediately yield 
the equality (15). This completes proof of Proposition 2.4.
   \end{proof}

\end{document}